\documentclass[twoside,11pt]{article}

\usepackage{jmlr2e}
\usepackage{url,verbatim}

\newcommand{\reals}{{\mbox{\bf R}}}
\newcommand{\eg}{{\it e.g.}}

\jmlrheading{17}{2016}{1-5}{8/15}{4/16}{Steven Diamond and Stephen Boyd}

\ShortHeadings{CVXPY: A Python-Embedded Modeling Language for Convex Optimization}{Diamond and Boyd}
\firstpageno{1}

\begin{document}

\title{CVXPY: A Python-Embedded Modeling Language for Convex Optimization}

\author{\name Steven Diamond \email diamond@cs.stanford.edu \\
       \name Stephen Boyd \email boyd@stanford.edu \\
       \addr Departments of Computer Science and Electrical Engineering\\
       Stanford University\\
       Stanford, CA 94305, USA}

\editor{Antti Honkela}

\maketitle

\begin{abstract}
CVXPY is a domain-specific language for convex optimization embedded in
Python.
It allows the user to express convex optimization problems in a natural
syntax that follows the math,
rather than in the restrictive standard form required by solvers.
CVXPY makes it easy to combine convex optimization with high-level
features of Python such as parallelism and object-oriented design.
CVXPY is available at \url{http://www.cvxpy.org/} under the GPL license,
along with documentation and examples.
\end{abstract}


\section{Introduction}

Convex optimization has many applications to fields as diverse as
machine learning, control, finance, and signal and
image processing \citep{BoV:04}.
Using convex optimization in an application requires either
developing a custom solver or converting
the problem into a standard form. Both of these tasks require
expertise, and are time-consuming and error prone.
An alternative is to use a domain-specific language (DSL)
for convex optimization, which allows the user to specify the
problem in a natural way that follows the math;
this specification is then automatically converted
into the standard form required by generic solvers.
CVX \citep{cvx}, YALMIP \citep{YALMIP},
QCML \citep{QCML}, PICOS \citep{PICOS}, and Convex.jl \citep{cvxjl}
are examples of such DSLs for convex optimization.

CVXPY is a new DSL for convex optimization.
It is based on CVX \citep{cvx},
but introduces new features such as signed disciplined
convex programming analysis and parameters.
CVXPY is an ordinary Python library,
which makes it easy to combine convex optimization with
high-level features of Python such as parallelism and
object-oriented design.

CVXPY has been downloaded by thousands of users and
used to teach multiple courses \citep{EE364a}.
Many tools have been built on top of CVXPY,
such as an extension for stochastic optimization \citep{DCSP}.

\section{CVXPY Syntax}

CVXPY has a simple, readable syntax inspired by CVX \citep{cvx}.
The following code constructs and solves a least squares problem
where the variable's entries are constrained to be between 0 and 1.
The problem data $A \in \reals^{m \times n}$ and $b \in \reals^m$
could be encoded as NumPy ndarrays or one of several other common
matrix representations in Python.

\begin{verbatim}
# Construct the problem.
x = Variable(n)
objective = Minimize(sum_squares(A*x - b))
constraints = [0 <= x, x <= 1]
prob = Problem(objective, constraints)

# The optimal objective is returned by prob.solve().
result = prob.solve()
# The optimal value for x is stored in x.value.
print x.value
\end{verbatim}

The variable, objective, and constraints
are each constructed separately and combined in the final problem.
In CVX, by contrast, these objects are created within the scope
of a particular problem.
Allowing variables and other objects to be created in isolation makes
it easier to write high-level code that
constructs problems (see \S\ref{OOCO}).

\section{Solvers}

CVXPY converts problems into a standard form known as conic form
\citep{NesNem:92}, a generalization of a linear program.
The conversion is done using graph implementations of convex
functions \citep{GB:08}.
The resulting cone program is equivalent to the original problem,
so by solving it we obtain a solution of the original problem.

Solvers that handle conic form are known as cone solvers; each one
can handle combinations of several types of cones.
CVXPY interfaces with the open-source cone solvers
CVXOPT \citep{CVXOPT},
ECOS \citep{bib:Domahidi2013ecos},
and SCS \citep{SCSpaper},
which are implemented in combinations of Python and C.
These solvers have different characteristics, such as the types of
cones they can handle and the type of algorithms employed.
CVXOPT and ECOS are interior-point solvers,
which reliably attain high accuracy for small and medium scale problems;
SCS is a first-order solver, which uses OpenMP to target multiple cores
and scales to large problems with modest accuracy.

\section{Signed DCP}

Like CVX, CVXPY uses disciplined convex programming (DCP) to verify
problem convexity \citep{GBY:06}.
In DCP, problems are constructed from a fixed library of functions
with known curvature and monotonicity properties.
Functions must be composed according to a simple set of rules
such that the composition's curvature is known.
For a visualization of the DCP rules, visit \url{dcp.stanford.edu}.

CVXPY extends the DCP rules used in CVX
by keeping track of the signs of expressions.
The monotonicity of many functions depends on the sign of their argument,
so keeping track of signs allows more compositions to be verified as convex.
For example, the composition \verb+square(square(x))+ would not be
verified as convex under standard DCP because the \verb+square+ function
is nonmonotonic.
But the composition is verified as convex under signed DCP
because \verb+square+ is increasing for nonnegative arguments
and \verb+square(x)+ is nonnegative.

\section{Parameters}

Another improvement in CVXPY is the introduction of parameters.
Parameters are constants whose symbolic properties
(\eg, dimensions and sign)
are fixed but whose numeric value can change.
A problem involving parameters can be solved repeatedly for different
values of the parameters without redoing computations that
do not depend on the parameter values.
Parameters are an old idea in DSLs for optimization,
appearing in AMPL \citep{AMPL}.

A common use case for parameters is computing a trade-off curve.
The following code constructs a LASSO problem \citep{BoV:04} where
the positive parameter $\gamma$ trades off the sum of squares error
and the regularization term.
The problem data are $A \in \reals^{m \times n}$ and $b \in \reals^m$.

\begin{verbatim}
x = Variable(n)
gamma = Parameter(sign="positive") # Must be positive due to DCP rules.
error = sum_squares(A*x - b)
regularization = norm(x, 1)
prob = Problem(Minimize(error + gamma*regularization))
\end{verbatim}

Computing a trade-off curve is trivially parallelizable,
since each problem can be solved independently.
CVXPY can be combined with Python multiprocessing
(or any other parallelism library) to
distribute the trade-off curve computation across many processes.

\begin{verbatim}
# Assign a value to gamma and find the optimal x.
def get_x(gamma_value):
    gamma.value = gamma_value
    result = prob.solve()
    return x.value

# Get a range of gamma values with NumPy.
gamma_vals = numpy.logspace(-4, 6)
# Do parallel computation with multiprocessing.
pool = multiprocessing.Pool(processes = N)
x_values = pool.map(get_x, gamma_vals)
\end{verbatim}

\section{Object-Oriented Convex Optimization}\label{OOCO}

CVXPY enables an object-oriented approach to constructing
optimization problems.
As an example, consider an optimal flow problem
on a directed graph $G = (V,E)$ with vertex set $V$ and (directed)
edge set $E$.
Each edge $e \in E$ carries a flow $f_e \in \reals$,
and each vertex $v \in V$ has an internal source that
generates $s_v \in \reals$ flow.
(Negative values correspond to flow in the opposite direction,
or a sink at a vertex.)
The (single commodity) flow problem is
(with variables $f_e$ and $s_v$)
\[
\begin{array}{ll} \mbox{minimize} & \sum_{e \in E}\phi_{e}(f_{e}) +
\sum_{v \in V}\psi_v(s_v), \\
\mbox{subject to} & s_v + \sum_{e \in I(v)}f_e = \sum_{e \in O(v)}f_e,
\quad \textrm{for all } v \in V,
\end{array}
\]
where the $\phi_e$ and $\psi_v$ are convex cost functions
and $I(v)$ and $O(v)$ give vertex $v$'s incoming and outgoing edges,
respectively.

To express the problem in CVXPY,
we construct vertex and edge objects, which store
local information such as optimization variables, constraints,
and an associated objective term.
These are exported as a CVXPY problem for each vertex and each edge.

\begin{verbatim}
class Vertex(object):
    def __init__(self, cost):
        self.source = Variable()
        self.cost = cost(self.source)
        self.edge_flows = []

    def prob(self):
        net_flow = sum(self.edge_flows) + self.source
        return Problem(Minimize(self.cost), [net_flow == 0])

class Edge(object):
    def __init__(self, cost):
        self.flow = Variable()
        self.cost = cost(self.flow)

    def connect(self, in_vertex, out_vertex):
        in_vertex.edge_flows.append(-self.flow)
        out_vertex.edge_flows.append(self.flow)

    def prob(self):
        return Problem(Minimize(self.cost))
\end{verbatim}

The vertex and edge objects are composed into a graph using the edges'
\verb+connect+ method.
To construct the single commodity flow problem,
we sum the vertices and edges' local problems.
(Addition of problems is overloaded in CVXPY
to add the objectives together and concatenate the constraints.)

\begin{verbatim}
prob = sum([object.prob() for object in vertices + edges])
prob.solve() # Solve the single commodity flow problem.
\end{verbatim}


\acks{We thank the many contributors to CVXPY.
This work was supported by DARPA XDATA.}


\newpage

\vskip 0.2in
\bibliography{15-408}

\end{document}